\input amssym.def
\input amssym
\magnification=1200
\parindent0pt
\hsize=16 true cm
\baselineskip=13  pt plus .2pt
$ $

\def\Z{\Bbb Z}
\def\A{\Bbb A}
\def\S{\Bbb S}
\def\R{\Bbb R}

\centerline {\bf On finite simple groups acting on homology spheres}

\bigskip  \bigskip

\centerline {Alessandra Guazzi*,  Bruno Zimmermann**}

\bigskip

\centerline {*SISSA} \centerline {Via Bonomea 256} \centerline {34136 Trieste, Italy}

\bigskip

\centerline {**Universit\`a degli Studi di Trieste} \centerline {Dipartimento di
Matematica e Informatica} \centerline {34127 Trieste, Italy}

\bigskip \bigskip

Abstract.  {\sl It is a consequence of the classical Jordan bound for finite subgroups
of linear groups that in each dimension $n$ there are only finitely many finite simple
groups which admit a faithful, linear action on the
$n$-sphere. In the present paper we prove an analogue for smooth actions on arbitrary
homology $n$-spheres: in each dimension  $n$ there are only finitely many finite simple
groups which admit a faithful, smooth action on some homology sphere of dimension $n$,
and in particular on the $n$-sphere.  We discuss also the finite simple groups which
admit an action on a homology sphere of dimension 3, 4 or 5. }

\bigskip

AMS-Classification-Numbers:   57S17;  57S25; 20E32; 20F65

\smallskip


\bigskip \bigskip

{\bf 1. Introduction}

\medskip

We consider finite groups, and in particular finite simple groups admitting
smooth actions on integer homology spheres. In the present paper, simple group
will always mean {\it nonabelian} simple group; also, all actions will be
smooth (or locally linear) and faithful.

\medskip

A motivation of the present paper is the basic problem of how to determine the finite
groups which admit an action on a sphere or homology sphere of a given dimension, and
the related problem to determine, for a given finite group, the minimal dimension of a
sphere or homology sphere on which it admits an action (comparing the situation with
that of linear actions on spheres, or real linear representations). For elementary
abelian $p$-groups, this is answered by classical Smith fixed point theory, and the
situation for arbitrary finite
$p$-groups is considered in [9]. At the other extreme, the question becomes
particularly attractive for finite simple groups. By [20, 14-16], the only finite
simple group which admits an action on a homology 3-sphere is the alternating group
$\A_5$, and the only finite simple groups acting on a homology 4-sphere are the
alternating groups $\A_5$ and $\A_6$; already these low-dimensional results require
heavy machinery from the theory of finite simple groups.

\medskip

The main result of the present  paper is the following:

\bigskip

{\bf Theorem 1.}  {\sl For each dimension $n$, up to isomorphism there are only finitely
many finite simple groups which admit an action on some homology sphere of dimension
$n$, and in particular on the $n$-sphere.}

\medskip

We note that any finite simple group admits many smooth actions on high-dimensional
spheres which are not equivalent to linear actions (see the survey [8, section 7]).

\medskip

The proof of Theorem 1 requires the full classification of the finite simple
groups. For the case of linear actions on spheres, Theorem 1 is a consequence
of the  existence of the classical Jordan number: for each dimension $n$ there
is an integer $j(n)$ such that each finite subgroup of ${\rm GL}_n(\Bbb C)$ has
a normal abelian subgroup of index at most $j(n)$. The optimal bound for all
$n$ has recently been determined in [6]: for $n \ge 71$ it is $(n+1)!$,
realized by the symmetric group $\S_{n+1}$ which is a subgroup of ${\rm
GL}_n(\Bbb C)$; this  requires again the classification of the finite simple
groups.  It is natural then to ask whether the Jordan number admits a
generalization for all finite groups acting on homology $n$-spheres, but this
remains open at present.

\medskip

The proof of Theorem 1 allows us to produce for each dimension $n$ a finite
list of finite simple groups which are the candidates for actions on homology
$n$-spheres; then one can identify those groups from the list which admit a linear action
on $S^n$ (or equivalently, have a faithful, real, linear representation in dimension
$n+1$), and try to  eliminate the remaining ones by refined methods.  We illustrate
this for dimension 5 by proving the following:

\bigskip

{\bf Theorem 2.}  {\sl  A finite simple group acting on a homology 5-sphere, and in
particular on the 5-sphere, is isomorphic to one of the following groups: an alternating
group $\A_5$, $\A_6$ or $\A_7$; a linear fractional group ${\rm PSL}_2(7)$; a unitary
group ${\rm  PSU}_4(2)$ or ${\rm PSU}_3(3)$.}

\bigskip

With the exception of the unitary group  ${\rm PSU}_3(3)$ these are exactly the finite
simple groups which admit a linear action on $S^5$. The alternating group
$\A_7$ acts on  $\R^7$ by permutation of coordinates, and hence also on a diagonal $\R^6$
and on $S^5$, and has the linear fractional group ${\rm PSL}_2(7)$ as a subgroup
(which contains the symmetric group $\S_4$ as a subgroup of index 7). The
unitary group  ${\rm PSU}_4(2)$, of order 25920, is a subgroup of index 2 in the Weyl
group of type $E_6$ which has an integer linear representation in dimension 6 (this is
one of the low-dimensional exceptions for the upper bound  $j(n) \le (n+1)!$ of the
Jordan number). The group ${\rm PSU}_3(3)$ admits a linear action on $S^6$ (see e.g.
[7]); we conjecture that it does not act on a homology 5-sphere but are not able to
exclude it at present.

\medskip

Note that, in contrast to Theorem 1 where one can neglect the first groups of each
infinite series in the classification, one has to exclude almost all finite simple
groups now; on the basis of Theorems 3 and 4 and Proposition 2 this is in general no
problem for the  large groups (e.g., the sporadic ones), but for some of the smaller
ones it may be more difficult (the worst case being that of ${\rm PSU}_3(3)$ which
remains open, in fact).

\medskip

Crucial for the proofs of Theorems 1 and 2 is also a control over the  minimal dimension
of an action of a linear fractional group ${\rm PSL}_2(p)$ and a linear group ${\rm
SL}_2(p)$:

\bigskip

{\bf Theorem 3.} {\sl For a prime $p \ge 5$, the minimal dimension of an action
of a linear fractional group ${\rm PSL}_2(p)$  on a mod $p$ homology sphere  is
$(p-1)/2$ if  $p  \equiv  1  \; {\rm mod} \; 4$, and $p-2$ if $p \equiv 3 \; {\rm mod}
\;  4$, and these are also lower bounds for the dimension of such an action of a linear
group  ${\rm SL}_2(p)$.}

\medskip

Whereas the groups ${\rm PSL}_2(p)$ admit linear actions on
spheres of the corresponding dimensions (see e.g. [10]), for the groups ${\rm
SL}_2(p)$ the minimal dimension of a linear action on a sphere is $p-2$ resp.
$p$ which is strictly larger than the lower bounds given in Theorem 2,
so the minimal dimension of an action on a homology sphere remains open here.

\medskip

Assuming Theorem 3, we will prove Theorems 1 and 2 in sections 2 and 3, and Theorem
3 in section 4.

\bigskip

{\bf 2. Proof of Theorem 1}

\medskip

In analogy with Theorem 3, we shall need the following result from Smith theory for
elementary abelian $p$-groups ([18]).

\bigskip

{\bf Theorem 4.}  {\sl  The minimal dimension of a faithful, orientation-preserving
action of an elementary abelian $p$-group $(\Bbb Z_p)^k$ on a mod $p$ homology sphere is
$k$ if $p = 2$, and $2k-1$ if $p$ is an odd prime.}

\bigskip

The {\it Proof of Theorem 1} is a consequence of Theorems 3 and 4 and the
classification of the finite simple groups.  Fixing a dimension $n$, we have to
exclude all but finitely many finite simple groups. By the classification of
the finite simple groups, a finite simple group is one of 26 sporadic groups,
or an alternating group, or a group of Lie type ([7, 2, 11, 12]).  We can
neglect the sporadic groups and have to exclude all but finitely many groups of
the infinite series.  Clearly, an alternating group $\A_m$  contains elementary
abelian subgroups $(\Bbb Z_2)^k$ of rank $k$ growing with the degree $m$, so
Theorem 4 excludes all but finitely many alternating groups and we are left
with the infinite series of groups of Lie type.  Since every group of Lie type
in characteristic $p \ge 5$ contains a subgroup isomorphic to either  ${\rm
SL}_2(p)$ or ${\rm PSL}_2(p)$, inside a root SL, the characteristic $p$ is
bounded by Theorem 3; also, root subgroups bound the field order by Theorem 4.
This excludes all but finitely many groups from each series of Lie type and
completes the proof of Theorem 1.

\bigskip

Remark.  The present short version of the proof of Theorem 1 was suggested by the
referee. In a preliminary version of the paper (arXiv:1106.1067v1), we went over the
list of the series of the finite simple groups of Lie type, excluding case by case all
but finitely many groups from each series. We thank the referee for his useful comments.

\bigskip

{\bf 3. Proof of Theorem 2}

\medskip

We shall present the proof of Theorem 2 for the most important case of the projective
linear groups  ${\rm PSL}_m(p^k)$; as in the proof of Theorem 1, for all other
classes of finite simple groups the proof reduces to this case, and we will
illustrate this only for the symplectic groups ${\rm PSp}_{2m}(p^k)$.  However,  we
present the  arguments for arbitrary dimension $n$ first, and then restrict to the case
$n = 5$ in a second step.

\medskip

With respect to the proof of Theorem 1, a basic additional tool for the proof of Theorem
2 is the Borel formula ([3, Theorem XIII.2.3]) which states that, for every elementary
abelian $p$-group $A \cong (\Z_p)^k$ acting on a mod $p$ homology $n$-sphere,
$$n-r(A)=\sum_H (r(H)-r(A));$$ the sum is taken over all subgroups $H \cong
(\Z_p)^{k-1}$ of index $p$ in $A$, and $r(H)$ and $r(A)$ are the dimensions of the fixed
point sets of $H$ and $A$, respectively (equal to  minus 1 if a fixed point set is
empty).

\bigskip

{\bf Lemma 1.}  {\sl Let $G$ be a finite group acting faithfuly on a homology
$n$-sphere, which has an elementary abelian $p$-subgroup $A \cong (\Z_p)^k$. Suppose
that one of the following conditions holds:

\smallskip

i) The fixed point sets of all subgroups $H$ of index $p$ in $A$ have the same dimension.

\smallskip

ii) All subgroups $H$ of index $p$ in $A$ are conjugate in $G$.

\smallskip

iii) All cyclic subgroups in $A$ are conjugate in $G$.

\smallskip

Then

\smallskip

\centerline  {$(p^k-1)/(p-1) \le  (n+1)/2$  \hskip 5mm if $p$ is odd,}

\medskip

\centerline {$2^k \le  n+2$  \hskip 5mm  if $p=2$.}}

\medskip

{\sl Proof.}  If all subgroups of index $p$ in $A$ fix submanifolds of the same
dimension, the Borel formula implies
$$n-r(A)=  {(p^k-1) \over (p-1)} (r(H)-r(A)).$$ We have  $r(A) \ge -1$, and $r(H) >
r(A)$ since the action of $G$ is faithful; moreover if $p$ is  odd, the fixed point set
of any $p$-group acting on a homology sphere is a mod $p$ homology sphere of even
codimension (since the action is orientation-preserving in this case), and hence $r(H) -
r(A)\geq 2$.  This proves case i) of the lemma, and then clearly also case ii). As for
iii), if all cyclic subgroups in $A$ are conjugate, their fixed point sets all have the
same dimension and hence, by an inductive use of the Borel formula, the same holds for
all subgroups $(\Z_p)^j$ of the same rank $j$ in $A$, so we again reduce to case i).
This completes the proof of the lemma.

\bigskip

Starting with the proof of Theorem 2 now, we distinguish  several cases.

\bigskip

{\bf 3.1.}  {\it Suppose that $G$ is a projective linear group  ${\rm PSL}_2(p^k)$  of
degree 2.}

\medskip

The subgroup of ${\rm PSL}_2(p^k)$ represented by triangular matrices with entries equal
to $1$ on the diagonal is isomorphic to an elementary abelian $p$-group $(\Z_p)^k$ and
is normalized by the cyclic group of diagonal matrices; in fact,

$$\pmatrix {
                 \omega^{-1}  &  0  \cr
                 0  &   \omega  \cr}
\pmatrix {
                 1  &   \sigma  \cr
                 0  &   1  \cr}
\pmatrix {
                 \omega  &  0  \cr
                 0  &   \omega^{-1}   \cr}  =
\pmatrix {
                 1  &   \omega^2\sigma  \cr
                 0  &   1  \cr},$$

for all $\sigma\in GF(p^k)$ (the finite Galois field with $p^k$ elements),
$\omega$ in $GF(p^k)^*$ (its multiplicative group).  In particular, no element
in $(\Z_p)^k$ is centralized by a diagonal matrix different from $\pm E^2$, and
the subgroup of upper triangular matrices is isomorphic to  $(\Z_p)^k\rtimes
\Z_{r}$, where $r=(p^k-1)/2$ if $p$ is odd, and $r = 2^k-1$ otherwise.

\medskip

If $p=2$ then clearly all $2^k-1$ involutions in $(\Z_2)^k$ are conjugate by elements in
$\Z_r = \Z_{2^k-1}$ and Lemma 1 implies that  $2^k \leq  n+2$. In particular if
$n=5$ we have $k\leq 2$, hence ${\rm PSL}_2(4)\cong \A_5$ is the only simple group of
type ${\rm PSL}_2(2^k)$ acting on a homology $5$-sphere.

\medskip

So we can assume that $p$ is an odd prime in the following.  Then exactly half
of the elements $\sigma$ of $GF(p^k)^* \cong \Z_{p^k-1}$ are squares,  and
accordingly there are exactly two conjugacy classes of elements of order $p$ in
$(\Z_p)^k < (\Z_p)^k\rtimes \Z_{r}$,  each with $(p^k-1)/2$ elements.

\medskip

Each cyclic subgroup $\Z_p$ of
$(\Z_p)^k$  is represented by $p$ matrices of the form

$$\pmatrix {
                 1  &   y \sigma  \cr
                 0  &   1  \cr},$$ for some fixed $\sigma\in GF(p^k)$, and with $y$
varying in the prime subfield $GF(p)$ of $GF(p^k)$.

\bigskip

If $k$ is odd, exactly half of the nontrivial elements of the prime field $GF(p)$ are
squares in $GF(p^k)$, so each cyclic subgroup $\Z_p$ contains elements of both conjugacy
classes of elements of order $p$; in particular, all such subgroups $\Z_p$ are conjugate
and hence Lemma 1 implies that $(p^k-1)/(p-1)  \le  (n+1)/2$.

\bigskip

If $k$ is even instead, all elements of the prime field $GF(p)$ are squares in
$GF(p^k)$, hence all nontrivial elements of each cyclic subgroup $\Z_p$ of
$(\Z_p)^k$ are conjugate, so there are exactly two conjugacy classes of the
$(p^k-1)/(p-1)$ cyclic subgroups $\Z_p$, each with the same number of elements.
In particular, for $k=2$ and $A=(\Z_p)^2 < {\rm PSL}_2(p^2)$, the Borel formula
implies
$$n-r(A) =  {p+1 \over 2}(r_1- r(A)) +  {p+1 \over 2}(r_2- r(A))$$ where $r_1$ and $r_2$
are the dimensions of the fixed point sets of the two conjugacy classes of
cyclic subgroups in $A$. We can suppose that $r_1 - r(A) \ne 0$; then $r_1 -
r(A) \ge 2$ since for odd $p$ actions are orientation-preserving and fixed
point sets have even codimensions. We conclude that  $p \le n$.

\medskip

Summarizing, for odd $p$ we have proved that

\bigskip

\centerline {$(p^k-1)/(p-1)  \le  (n+1)/2$  \hskip 5mm if $k$ is odd,}

\medskip

\centerline {$p \le n$ \hskip 5mm  if $k$ is even.}

\bigskip

Specializing to $n=5$ now, Theorems 3 and 4 imply that  $p
\leq  7$ and $k \leq 3$.  So for $k = 1$ we are left with the groups ${\rm PSL}_2(5)$
and ${\rm PSL}_2(7)$ of Theorem 2 which admit linear actions on $S^5$.  The case $k = 3$
is excluded by the first of the preceding inequalities, and
$k = 2$ leaves only the possibilities $p=3$ and $p=5$. The group
${\rm PSL}_2(9) \cong \A_6$ is one of the groups of Theorem 2 admitting a linear action
on $S^5$, so it remains to exclude ${\rm PSL}_2(25)$.

\medskip

Suppose that ${\rm PSL}_2(25)$ acts on a homology 5-sphere. The Borel formula implies
easily that the two conjugacy classes of cyclic subgroups of $(\Z_5)^2$ have empty resp.
$1$-dimensional fixed point sets (the case $r_1 = -1$ and $r_2 = 1$). Consider the
subgroup  $C \cong (\Z_5 \times
\Z_5)\rtimes \Z_4$ of the upper triangular matrices, with a faithful action of $\Z_4$ on
each cyclic subgroup
$\Z_5$. Let  $H\cong
\Z_5$ be a subgroup of $C$ with 1-dimensional fixed point set, homeomorphic to the
1-sphere $S^1$. Since
$H$ is normal in $C$, also $C$ acts on $S^1$,  with some kernel $K \geq H$ not
containing the whole subgroup $\Z_5\times \Z_5$. Since $K/H$ is a normal subgroup of $C/H
\cong
\Z_5 \rtimes \Z_4$, it follows that $K=H$, hence $\Z_5\rtimes \Z_4$ acts faithfully on
$S^1$ which is a contradiction.

\medskip

This excludes the group ${\rm PSL}_2(25)$ and completes the proof of Theorem 2 for the
groups ${\rm PSL}_2(p^k)$.

\bigskip

{\bf 3.2.} {\it Suppose now that $G$ is a projective linear group  ${\rm PSL}_m(p^k)$ of
degree $m \geq 3$.}

\medskip

The group ${\rm PSL}_m(p^k)$ admits an elementary abelian
$p$-subgroup $P \cong (\Z_p)^{(m-1)k}$, represented by all upper triangular matrices of
the form
$$M(v) =   \pmatrix {
                 E_{n-1}  &   v    \cr
                 0        &   1            \cr}$$

where $v$ is a (column-) vector in the vector space $(GF(p^k))^{m-1}$ (noting that
$M(v) \cdot M(w) = M(v + w)$).  Theorem 4 implies now that

\bigskip

\centerline {$(m-1)k\leq (n+1)/2$  \hskip 5mm  if $p$ is odd,}

\medskip

\centerline {$m(k-1)\leq n$  \hskip 5mm if $p=2$.}

\bigskip

Considering the linear group  ${\rm SL}_{m-1}(p^k)$ as a subgroup of ${\rm PSL}_m(p^k)$
in the standard way, this normalizes $P$ since
$$\pmatrix {
                 A^{-1}   &   0    \cr
                 0        &   1   \cr} M(v)
\pmatrix {
                 A  &    0    \cr
                 0  &    1   \cr}  =  M(A^{-1}v),$$ for each $A \in {\rm
SL}_{m-1}(p^k)$.  Hence ${\rm PSL}_m(p^k)$ has a subgroup isomorphic to the semidirect
product $P \rtimes {\rm SL}_{m-1}(p^k) \cong {\rm ASL}_{m-1}(p^k)$, the group of affine
transformations in dimension $n-1$ over the field  $GF(p^k)$. Considering $P$ as a
vector space over the field
$GF(p^k)$,  under the action of ${\rm SL}_{m-1}(p^k)$ every two linear subspaces of the
same dimension are conjugate.  In particular, if $k=1$ every two subgroups of the same
index in $P$ are conjugate, so Lemma 1 implies that

\bigskip

\centerline {$(p^{m-1}-1)/(p-1) \leq  (n+1)/2$  \hskip 5mm  if $p$ is odd,}

\medskip

\centerline {$2^{m-1}\leq n+2$  \hskip 5mm  if $p=2$.}

\bigskip

Suppose now that $n=5$.  If $p$ is odd we have $p+1 \leq (p^{m-1}-1)/(p-1) \leq 3$ which
has no solution in $p$.  For $p=2$ instead it follows that $2^{m-1}\leq 7$, hence $m=3$
and, since $2k=(m-1)k\leq 5$, we have $k=1$ or 2. The group ${\rm PSL}_3(2) \cong  {\rm
PSL}_2(7)$ is one of the groups in Theorem 2, acting linearly on
$S^5$. The remaining possibility is ${\rm PSL}_3(4)$ which has one conjugacy class of
involutions (see e.g. [7]); in particular, all involutions in its elementary abelian
$2$-subgroup $(\Z_2)^4$ are conjugate and hence by Lemma 1, ${\rm PSL}_3(4)$ does not act
on a homology $5$-sphere.

\medskip

This completes the proof of Theorem 2 for the projective linear groups
${\rm PSL}_m(p^k)$.

\medskip

Incidentally, since $\A_8 \cong {\rm PSL}_4(2)$, this shows also that the alternating
group $\A_m$ does not act on a homology 5-sphere if $m \ge 8$.

\medskip

Similar as the proof of Theorem 1, also the proof of Theorem 2 is based on the case of
the projective linear groups and the linear groups.  Since Theorem 3 is most likely not
best possible for the linear groups ${\rm SL}_2(p)$, the following improvement for
dimension 5 turns out to  be useful.

\bigskip

{\bf Proposition 1.}  {\sl  If a linear group  ${\rm SL}_2(q)$ acts faithfully on a
homology $5$-sphere $M^5$ then $q\leq 5$.}

\medskip

{\it Proof.} If $q$ is a power of $2$ then ${\rm SL}_2(q) = {\rm PSL}_2(q)$ and the
claim  has already been proved in 3.1, so we can assume that $q$ is odd; then
${\rm SL}_2(q)/Z \cong {\rm PSL}_2(q)$ where $Z \cong \Z_2$ denotes the center of
${\rm SL}_2(q)$ (consisting of the matrices $\pm E_2$).  The $2$-Sylow subgroups of ${\rm
SL}_2(q)$  are (generalized)  quaternion groups, and consequently ${\rm SL}_2(q)$
contains a copy of the quaternion group $Q_8$ of order 8, with center $Z$. Since a
quaternion group has periodic cohomology of period 4, it cannot act freely on a
homology sphere of dimension $4m+1$ (see [4, section VI.9]), and in particular not on a
homology 5-sphere. So there is a nontrivial element in $Q_{8}$ with nonempty fixed point
set, and hence also the central involution in $Q_8 < {\rm SL}_2(q)$ has nonempty fixed
point set. This fixed point set is  a mod 2 homology sphere of even codimension, i.e.
either a mod 2 homology 3-sphere or a 1-sphere $S^1$.

\medskip

Suppose now that $Z$ has fixed point set $S^1$, and let $K$ be the normal subgroup of
${\rm SL}_2(q)$ fixing  $S^1$ pointwise; then either $K = Z$ or $K = {\rm SL}_2(q)$. In
the first case, the factor group
${\rm SL}_2(q)/K \cong  {\rm PSL}_2(q)$ acts faithfully on $S^1$ which is  possible only
for $q = 3$. In the second case, ${\rm SL}_2(q)$ acts orthogonally on a tubular
neighbourhood of $S^1$ in $M^ 5$ and hence on $S^3$ (see [4, section III]) which is
possible only for $q=5$.

\medskip

Suppose then that the fixed point set of $Z$ is a  mod 2 homology
$3$-sphere $M^3$; then ${\rm SL}_2(q)/Z \cong {\rm PSL}_2(q)$  and its subgroup
$Q_8/Z \cong (\Z_2)^2$ act faithfully on $M^3$. Since ${\rm PSL}_2(q)$ has only one
conjugacy class of involutions, the Borel formula  implies that the fixed point set of
$(\Z_2)^2$ on $M^3$ is a 0-sphere $S^0$, i.e. two points which are fixed points also of
the action of $Q_8$ on $M^5$. Hence
$Q_8$ acts linearly on the border $S^4$ of a regular neighbourhood in $M^5$ of either
one of these two fixed points, and the central involution fixes a 2-sphere $S^2 = M^3
\cap S^4$ invariant under the action of $Q_8$. Then $Q_8$ acts linearly also on the
orthogonal complement $S^1$ of this 2-sphere in $S^4$, and hence is a subgroup of ${\rm
O}(3) \times {\rm O}(2)$. However it is easy to check that $Q_8$ is not a subgroup of
${\rm O}(3) \times {\rm O}(2)$, so we get a contradiction and this case does not occur.

\medskip

This completes the proof of Proposition 1.

\medskip

Finally, as an application of the previous methods we consider the case of the symplectic
groups:

\bigskip

{\bf 3.3.}  {\it  Suppose that $G$ is a  symplectic group  ${\rm PSp}_{2m}(p^k)$.}

\medskip

If $A$ is any matrix in the linear group ${\rm SL}_m(p^k)$, then

$$M(A)=  \pmatrix {
                 A  &      0          \cr
                 0  &    ^tA^{-1}   \cr}$$

is a matrix in the symplectic group ${\rm Sp}_{2m}(p^k)$ (leaving invariant the standard
symplectic form). Hence ${\rm SL}_m(p^k)$ injects into
${\rm Sp}_{2m}(p^k)$, and the projection ${\rm Sp}_{2m}(p^k) \to {\rm PSp}_{2m}(p^k)$
has no kernel when restricted to the affine group ${\rm ASL}_{m-1}(p) < {\rm SL}_m(p^k)$
(see 3.2).

\medskip

Suppose that ${\rm PSp}_{2m}(p^k)$ acts on a homology 5-sphere. If $p$ is odd,
we have proved in section 3.2 that $(p^{m-1}-1)/(p-1)\leq 3$, so $m=2$. Also, since
${\rm PSp}_4(q) > {\rm Sp}_2(q) \cong {\rm SL}_2(q)$, Proposition 1 implies
that $q \le 5$.  It is then enough to note that ${\rm PSp}_4(3) \cong {\rm
PSU}_4(2)$ acts linearly on $S^5$, while ${\rm PSp}_4(5) > {\rm PSL}_2(25)$
does not act on a homology 5-sphere.

\medskip

If instead $q=2^k$ then $2^{m-1}\leq 7$ and
$(m-1)k\leq 5$, so either $m=2$ and $k=1,2$, or $m=3$ and $k=1$. Hence it is enough to
note that ${\rm PSp}_4(2)$ is not simple, and that ${\rm PSp}_4(4) > {\rm PSL}_2(16)$
and ${\rm PSp}_6(2) > A_8 \cong {\rm PSL}_4(2)$ (see [7]) do not act on a homology
5-sphere.

\medskip

This completes the proof of Theorem 2 for the symplectic groups; for all other classes
of simple groups, similar methods apply so we will not present the details here.

\bigskip

{\bf 4.  Proof of Theorem 3}

\medskip

In the following, for an odd prime $p$ and a positive integer $q$, we denote by
$\Z_p \rtimes \Z_q$ a semidirect product or split metacyclic group, with normal subgroup
$\Z_p$ on which $\Z_q$ acts by conjugation; we assume that a generator of $\Z_q$ acts on
$\Z_p$ by multiplication with $y \in  \Z_p^*$.

\medskip

Suppose that $\Z_p \rtimes \Z_q$ acts faithfully and orientation-preservingly  on a mod
$p$ homology  $n$-sphere, that is  on a closed orientable $n$-manifold $M$ with the mod
$p$ homology of the
$n$-sphere. We denote by $F$ the fixed point set of $\Z_p$, invariant under the action
of  $\Z_p \rtimes \Z_q$; by Smith fixed point theory, $F$ is a mod $p$ homology
$d$-sphere of even codimension
$n-d$ (since the action of $\Z_p$ is orientation-preserving; as usual $d = -1$ if $F$ is
empty). Let $M_0 = M - F$ denote the complement of $F$ in $M$; then $M_0$ is an open
$n$-manifold which, by Alexander-Lefschetz  duality, has the mod $p$ cohomology of a
sphere of odd  dimension $m = n - d - 1$ (see e.g. [17, section 71]). Also, $\Z_p
\rtimes \Z_q$ acts on $M_0$ with $\Z_p$ acting freely, so
$M_0 \to M_0/\Z_p$ is a regular covering with covering group $\Z_p$, and the action of
$\Z_q$ on $M_0$ (not necessarily orientation-preserving) projects to an action of
$\Z_q$ on  $M_0/\Z_p$.

\bigskip

{\bf Lemma 2.} {\sl   With notations as above,  $y^{(m+1)/2} =  \pm 1$.}

\medskip

{\it Proof.}  We apply the cohomology spectral sequence associated to the regular
$\Z_p$-covering  $M_0 \to M_0/\Z_p$ (see [13, 5, 4]);  this takes the form
$$E_2^{i,j}  =  H^i(\Z_p; H^j(M_0;K))   \;\;\;  \Rightarrow  \;\;\;
H^{i+j}(M_0/\Z_p;K),$$ i.e. it converges to the graded group associated to a filtration
of
$H^*(M_0/\Z_p;K)$; here $K$ denotes an abelian coefficient group or a commutative ring.

\medskip

We consider first integer coefficients, in general without indicating them. Since
$M_0$ has the mod $p$ homology and cohomology of a sphere of dimension $m = n - d - 1$,
in dimensions $j \ne 0, m$ the groups $H^j(M_0)$ are torsion of order coprime to
$p$, and  $H^m(M_0) \cong H^0(M_0) \cong \Z$. As a consequence,
$H^i(\Z_p; H^j(M_0))$ is trivial for $j \ne 0, m$, and
$H^i(\Z_p; H^m(M_0)) \; \cong \; H^i(\Z_p; H^0(M_0)) \; \cong \; H^i(\Z_p; \Z)$   is
isomorphic to $\Z$ for $i=0$, to $\Z_p$ if $i>0$ is even, and trivial if $i$ is odd (see
[4] or [13] for basic facts about the cohomology of groups). Hence the spectral
sequence is concentrated in the rows $i=0$ and $m$, and the only possibly nontrivial
differentials are
$d_{m+1}^{i,m}: E_{m+1}^{i,m} \to E_{m+1}^{i+m+1,0}$, of bidegree $(m+1,-m)$, starting
with
$$d_{m+1}^{0,m}: H^0(\Z_p; H^m(M_0)) \cong \Z  \;\;  \to   \;\;     H^{m+1}(\Z_p;
H^0(M_0)) \cong \Z_p.$$  Passing to the limit
$H^{i+j}(M_0/\Z_p)$ this easily implies that $H^m(M_0/\Z_p) \cong \Z$, and
$\Z_q$ acts by multiplication with $\pm 1$ on the latter group (multiplication with +1
if $\Z_q$ acts orientation-preservingly on $M_0$ and hence on $M_0/\Z_p$).

\bigskip

Next we consider the spectral sequence with coefficients in $\Z_p$; now
$H^j(M_0;\Z_p)$ is isomorphic to $\Z_p$ for $j=m$ and 0, and trivial otherwise. Again
the spectral sequence is concentrated in the rows  $j=m$ and 0 and, for all $i$,
$$H^i(\Z_p; H^m(M_p;\Z_p)) \; \cong \; H^i(\Z_p; H^0(M_0;\Z_p)) \; \cong \; H^i(\Z_p;
\Z_p) \; \cong \Z_p.$$

\medskip

In fact, the cohomology ring $H^*(\Z_p;\Z_p)$ is the tensor product of a polynomial
algebra $\Z_p[x_2]$ on a 2-dimensional generator $x_2$ and an exterior algebra
$E(x_1)$ on a 1-dimensional generator $x_1$ (see [1, Corollary II.4.2]); also, $x_2$ is
the image of $x_1$ under the mod $p$ Bockstein homomorphism.

\medskip

Since $H^k(M_0/\Z_p;Z_p)$ is trivial if $k$ is larger than the dimension $n$ of
$M_0/\Z_p$, the differentials  $d_{m+1}^{i,m}$ are isomorphisms if $i+j>n$, and this
implies that all differentials $d_{m+1}^{i,m}$ are isomorphisms, starting with
$$d_{m+1}^{0,m}: H^0(\Z_p; H^m(M_0;\Z_p)) \cong \Z_p \to  H^{m+1}(\Z_p; H^0(M_0;\Z_p))
\cong \Z_p.$$ In fact, for any coefficients $K$, the cyclic group
$\Z_p$ has periodic (Tate) cohomology $\hat H^*(\Z_p;K)$ of period 2, the isomorphisms
being given by the cup product with an element $u \in \hat H^2(\Z_p;H^0(M_0)) \cong
H^2(\Z_p;\Z)$ (see [4, chapter IV.9]); this periodicity applies to the rows $j=0$ and
$j=m$ of the spectral sequence. Since the differentials
$d_{m+1}^{i,m}$ are isomorphisms for large $i$, shifting dimensions by the cup product
with $u$ shows that they are isomorphisms for all
$i \ge 0$, by the multiplicative structure of spectral sequences, with cup-product
pairing
$$H^i(\Z_p; H^j(M_0;K)) \otimes H^{i'}(\Z_p; H^{j'}(M_0;K')) \to  H^{i+i'}(\Z_p;
H^{j+j'}(M_0;K \otimes K'))$$  (see [19, section 9.4] for the multiplicative properties
of the spectral sequence associated to a fibration; note that the spectral sequence of a
covering can be obtained form that of a suitable fibration, see [4], chapter VII.7 and
in particular exercise 3 for the homology version).

\medskip

Passing to the limit, it follows that
$H^m(M_0/\Z_p;Z_p) \cong H^m(\Z_p;H^0(M_0;\Z_p)) \cong \Z_p$. Since $\Z_q$ acts by
multiplication with $\pm 1$ on
$H^m(M_0/\Z_p;Z) \cong \Z$, it acts by multiplication with
$\pm 1$ also on $H^m(M_0/\Z_p;Z_p) \cong \Z_p$.

\medskip

A generator of $\Z_q$ acts by multiplication with $y \in  \Z_p^*$ also on
$H_1(\Z_p) \cong \Z_p$, and in the same way on $H^1(\Z_p;H^0(M_0;\Z_p))  \cong
H^1(\Z_p;\Z_p)
\cong
\Z_p$ (generated by $x_1$); via the Bockstein homomorphism, it has the same action on
$H^2(\Z_p;H^0(M_0;\Z_p))  \cong H^2(\Z_p;\Z_p)
\cong \Z_p$ (generated by $x_2$). Then it acts on a generator $x_1x_2^{(m-1)/2}$ of
$H^m(\Z_p; H^0(M_0;\Z_p))  \cong H^m(\Z_p; \Z_p) \cong \Z_p$ by multiplication with
$y^{(m+1)/2}$, and hence by multiplication with
$y^{(m+1)/2}$ also on  $H^m(M_0/\Z_p; Z_p) \cong H^m(\Z_p;H^0(M_0;\Z_p))$. On the other
hand, it acts on $H^m(M_0/\Z_p; Z_p)$ by multiplication with $\pm 1$, hence
$y^{(m+1)/2} =  \pm 1$.

\medskip

This completes the proof of Lemma 2.

\bigskip

{\it Proof of Theorem 2.}  We consider the case of ${\rm SL}_2(p)$ first. The subgroup
of ${\rm SL}_2(p)$ of all upper triangular matrices is a semidirect product $\Z_p
\rtimes \Z_q$, with $q=p-1$ ($\Z_q$ is the subgroup of all diagonal matrices, $\Z_p$ of
all triangular matrices with both diagonal entries equal to 1). The subgroup $\Z_q$
contains the central subgroup $\pm E_2$ of
${\rm SL}_2(p)$, and all other elements of $\Z_q$ act nontrivially on $\Z_p$. Hence, with
notations as above, $y$ has order $q/2 = (p-1)/2$.

\medskip

By Lemma 2, $y^{(m+1)/2} =  \pm 1$, with $m = n-d-1$. Suppose first that $\Z_p$ has
empty fixed point set, so $M_0 = M$, $d=-1$ and
$m=n$.  Since $\Z_q$ acts orientation-preservingly on $M$, $y^{(n+1)/2} = 1$ (by Lemma 2
and its proof). Now the order $q/2$ of $y$ divides $(n+1)/2$ which implies $n
\ge q-1 = p-2$, proving Theorem 2 in the present case.

\medskip

Now suppose that $\Z_p$ does not act freely, so $F$ has dimension $d \ge 0$. If
$y^{(m+1)/2} = 1$ then $q/2$ divides  $(m+1)/2$, hence $m+1 = n-d \ge q = p-1$ and $n
\ge p-1$. On the other hand, if $y^{(m+1)/2} = -1$ then
$y^{m+1} = 1$, so $q/2$ divides $m+1 = n-d$ which implies $n \ge q/2 = (p-1)/2$; note
that in this case the order $q/2 = (p-1)/2$ of $y$ has to be even, so $p \;  \equiv
\;  1  \; {\rm mod} \;  4$.  This finishes the proof of Theorem 2 in the case of
${\rm SL}_2(p)$.

\medskip

The group ${\rm PSL}_2(p) = {\rm SL}_2(p)/\pm E_2$ has a subgroup $\Z_p \rtimes
\Z_q$, with  $q=(p-1)/2$, and $\Z_q$ acts effectively on $\Z_p$. The order
$q=(p-1)/2$ of $y$ is the same now as in the previous case and the same proof applies,
completing the proof of Theorem 2.

\bigskip

In a similar way one obtains:

\bigskip

{\bf Proposition 2.} {\sl  For a prime $p \ge 5$ and a positive integer $q$,
let $\Z_p \rtimes \Z_q$ denote a semidirect product with an effective action of
$\Z_q$ on the normal subgroup $\Z_p$. The minimal dimension of a smooth,
orientation-preserving action of $\Z_p \rtimes \Z_q$ on a mod $p$ homology
$n$-sphere is $2q-1$ if  $q$ is odd, and $q$ if $q$ is even.}

\bigskip

As in the case of ${\rm PSL}_2(p)$, the values in Proposition 2 coincide with
the minimal dimension of a linear, orientation-preserving action on a sphere;
note also that the lower bounds in Proposition 2 determine those for the groups
${\rm PSL}_2(p)$ in Theorem 2, with $q=(p-1)/2$.

\bigskip

\centerline {\bf References}

\bigskip

\item {[1]}  A. Adem, R.J. Milgram, {\it Cohomology of finite groups.} Grundlehren der
math. Wissenschaften 309, Springer-Verlag 1994

\item {[2]} M. Aschbacher, {\it Finite Group Theory.} Cambridge University Press 1986

\item {[3]} A. Borel, {\it Seminar on Transformation Groups.} Annals of Math. Studies
46, Princeton University Press 1960

\item {[4]} K.S. Brown, {\it Cohomology of Groups.}  Graduate Texts in Mathematics 87,
Springer 1982

\item {[5]} H. Cartan, S. Eilenberg, {\it Homological Algebra.} Princeton University
Press, Princeton 1956

\item {[6]}  M.J. Collins, {\it  On Jordan's theorem for complex linear groups.}
J. Group Theory 10,  411-423 (2007)

\item {[7]} J.H. Conway, R.T. Curtis, S.P. Norton, R.A. Parker, R.A.Wilson, {\it Atlas
of Finite Groups.} Oxford University Press 1985

\item {[8]}  M.W. Davis, {\it  A survey of results in higher dimensions.}  The Smith
Conjecture, edited by J.W. Morgan, H. Bass, Academic Press 1984,  227-240

\item {[9]}  R.M. Dotzel, G.C. Hamrick, {\it  $p$-group actions on homology spheres.}
Invent. math. 62, 437-442  (1981)

\item {[10]} W. Fulton, J. Harris, {\it Representation Theory: A First Course.} Graduate
Texts in Mathematics 129,  Springer-Verlag 1991

\item {[11]} D. Gorenstein, {\it The Classification of Finite Simple Groups.} Plenum
Press, New York 1983

\item {[12]} D. Gorenstein, {\it Finite Simple Groups: An Introduction to  their
Classification.}  Plenum Press, New York 1982

\item {[13]} S. MacLane, {\it Homology.} Springer-Verlag,  Berlin 1963

\item {[14]} M. Mecchia, B. Zimmermann, {\it On finite groups acting on
$\Z_2$-homology 3-spheres.} Math. Z. 248, 675-693 (2004)

\item {[15]} M. Mecchia, B. Zimmermann, {\it On finite simple groups acting on integer
and mod 2 homology 3-spheres.}  J. Algebra 298, 460-467  (2006)

\item {[16]} M. Mecchia, B. Zimmermann, {\it On finite simple and nonsolvable groups
acting on homology 4-spheres.} Top. Appl. 153,  2933-2942  (2006)

\item {[17]} J.R. Munkres, {\it Elements of Algebraic Topology.} Addison-Wesley
Publishing Company 1984

\item {[18]}  P.A. Smith, {\it Permutable periodic transformations.} Proc. Nat. Acad.
Sci. U.S.A. 30, 105 - 108 (1944)

\item {[19]} E.H. Spanier, {\it Algebraic Topology.}  McGraw-Hill, New York 1966

\item {[20]} B. Zimmermann, {\it On finite simple groups acting on homology 3-spheres.}
Top. Appl. 125, 199-202 (2002)

\bye